\documentclass[12pt]{article}
\usepackage{amsmath, amssymb, amsfonts, amsthm}
\def\CC{{\rm \kern.24em \vrule width.02em height1.4ex
depth-.05ex \kern-.26em C}}

\def\TagOnRight

\def\AA{{\it I}\hskip-3pt{\tt A}}

\def\QQ{\rlap {\raise 0.4ex \hbox{$\scriptscriptstyle |$}}
  {\hskip -0.1em Q}}

\newcommand{\be}{\begin{equation}}
\newcommand{\ee}{\end{equation}}
\newcommand{\bea}{\begin{eqnarray}}
\newcommand{\eea}{\end{eqnarray}}
\newcommand{\Bea}{\begin{eqnarray*}}
\newcommand{\Eea}{\end{eqnarray*}}

\newcommand{\bi}{\begin{itemize}}
\newcommand{\ei}{\end{itemize}}

\newtheorem{Definition}{Definition}
\newtheorem{Theorem}[Definition]{Theorem}
\newtheorem{Lemma}[Definition]{Lemma}

\newtheorem{Corollary}[Definition]{Corollary}

\newtheorem{Revised Question}[Definition]{Revised Question}

\theoremstyle{remark}

\begin{document}

\title{A tightness criterion for homology manifolds with or without boundary}
\author{Bhaskar Bagchi\\
Theoretical Statistics and Mathematics Unit\\
Indian Statistical Institute\\
Bangalore - 560 059, India.\\
email: bbagchi@isibang.ac.in}
\date{}

\maketitle

\abstract{A simplicial complex $X$ is said to be tight with respect to a field $\mathbb{F}$ if $X$ is connected and, for every induced subcomplex $Y$ of $X$, the linear map $H_\ast (Y; \mathbb{F}) \rightarrow H_\ast (X; \mathbb{F})$ (induced by the inclusion map) is injective. This notion was introduced by K\"{u}hnel in [10]. In this paper we prove the following two combinatorial criteria for tightness.
(a) Any $(k+1)$-neighbourly $k$-stacked $\mathbb{F}$-homology manifold with boundary is $\mathbb{F}$-tight. Also, (b) any $\mathbb{F}$-orientable $(k+1)$-neighbourly $k$-stacked $\mathbb{F}$-homology manifold without boundary is $\mathbb{F}$-tight, at least if its dimension is not equal to $2k+1$.

The result (a) appears to be the first criterion to be found for tightness of (homology) manifolds with boundary. Since every $(k+1)$-neighbourly $k$-stacked manifold without boundary is, by definition, the boundary of a $(k+1)$-neighbourly $k$-stacked manifold with boundary - and since we now know several examples (including two infinite families) of triangulations from the former class - theorem (a) provides us with many examples of tight triangulated manifolds with boundary.

The second result (b) generalizes a similar result from [2] which was proved for a class of combinatorial manifolds without boundary. We believe that theorem (b) is valid for dimension $2k+1$ as well. Except for this lacuna, this result answers a recent question of Effenberger [8] affirmatively.}

We recall that, for any simplicial complex $X$ of dimension $d$, and for $0 \leq k \leq d$, the $k$-skeleton $\text{Skel}_k(X)$ of $X$ is the subcomplex $\{\alpha \in X: \dim (\alpha) \leq k\}$. An $\mathbb{F}$-homology ball $B$ of dimension $d+1$ is said to be $k$-stacked if $\text{Skel}_{d-k}(B)=\text{Skel}_{d-k}(\partial B)$. An $\mathbb{F}$-homology sphere of dimension $d$ is said to be $k$-stacked if there is a $k$-stacked $\mathbb{F}$-homology ball $B$ (of dimension $d+1$) such that $S=\partial B$.

This definition was extended in [13] to homology manifolds with or without boundary in the obvious way. Thus, an $\mathbb{F}$-homology manifold $\Delta$ with boundary, say of dimension $d+1$, is $k$-stacked if $\text{Skel}_{d-k}(\Delta)=\text{Skel}_{d-k}(\partial \Delta)$. An $\mathbb{F}$-homology manifold without boundary is $k$-stacked if it is the boundary of an $\mathbb{F}$-homology $k$-stacked manifold with boundary.

Throughout the rest of this paper, $k$ is a strictly positive integer. For any simplicial complex $X$, its vertex set will be denoted by $V(X)$. For any subset $A$ of $V(X)$, $X[A]:= \{\alpha \in X: \alpha \subseteq A\}$ is the induced subcomplex of $X$ with vertex set $A$. For $0 \leq i \leq d$, $f_i (X) := \# \{\alpha \in X: \dim (\alpha)=i\}$. In particular, $f_0(X)=\# V(X)$. Also, $\beta_i (X; \mathbb{F})=\dim_{\mathbb{F}} H_i (X; \mathbb{F})$ and $\tilde \beta_i (X; \mathbb{F})=\dim_{\mathbb{F}} \tilde H_i (X; \mathbb{F})$ are the Betti numbers and the reduced Betti numbers (respectively) of $X$. Thus $\beta_i =\tilde \beta_i +\delta_{i0}$.

The sigma vector $(\sigma_0, \ldots,\sigma_d)$ of a $d$-dimensional simplicial complex $X$ (w.r.t. $\mathbb{F}$) was defined in [2] by the formula
\begin{eqnarray*}
\sigma_i = \sigma_i (X; \mathbb{F}) := \sum\limits_{A \subseteq V(X)} \frac{\tilde \beta_i (X[A]; \mathbb{F})}{\binom{f_0(X)}{\# (A)}}, ~0 \leq i \leq d.
\end{eqnarray*}
We also put $\sigma_i (X; \mathbb{F})=0$ if $i <0$ or $i >d$.

The new ingredient in this paper is the following lemma. Note that its proof closely follows the proof of Lemma 2.11 in [1], which in turn paraphrases the proof of Theorem 2.3 (ii) in [12].

\begin{Lemma} If $B$ is a $k$-stacked $\mathbb{F}$-homology ball then $\sigma_i (B; \mathbb{F})=0$ for $i \geq k \geq 1$.
\end{Lemma}

{\bf Proof:} Let $\dim (B) =d+1$ and $S=\partial B$. Thus $S$ is an $\mathbb{F}$-homology $d$-sphere with $\text{Skel}_{d-k}(S)=\text{Skel}_{d-k}(B)$. Take a new vertex $x \not\in V := V(B)$ and put $\hat B=x \ast B$, the cone over $B$ from the vertex $x$. Let $\hat S=\partial \hat B=B\cup (x \ast S)$. Let $\hat V=V \sqcup \{x\}=V(\hat B)$. Thus $\hat B$ is an $\mathbb{F}$-homology $(d+2)$-ball, and $\hat S$ is an $\mathbb{F}$-homology $(d+1)$-sphere.

Take any subset $\alpha$ of $V$. Put $\beta =\hat V \backslash \alpha, \gamma =V\backslash \alpha$. Thus $\beta = \gamma \sqcup \{x\}$. Since $d-i \leq d-k$, we have $\text{Skel}_{d-i}(B)=\text{Skel}_{d-i}(S)$, and hence $\text{Skel}_{d-i}(\hat S[\beta])=\text{Skel}_{d-i} (x\ast S[\gamma])$. Since $x\ast S[\gamma]$ is homologically trivial and $x\ast S \subseteq \hat S$, this implies (in the usual notation of homology theory) that we have
\begin{eqnarray*}
Z_{d-i} (x\ast S[\gamma]) &=& B_{d-i} (x\ast S[\gamma]) \subseteq B_{d-i} (\hat S[\beta])\\
\subseteq Z_{d-i} (\hat S[\beta]) &=& Z_{d-i} (x\ast S[\gamma]).
\end{eqnarray*}
Therefore, we get $B_{d-i} (\hat S[\beta])=Z_{d-i} (\hat S[\beta])$, and hence $H_{d-i}(\hat S[\beta])=0$. Since $\hat S$ is an $\mathbb{F}$-homology $(d+1)$-sphere, and $\beta$ is the complement of $\alpha$ in the vertex set of $\hat S$, Alexander duality and the exact sequence for pairs imply that $H_i (\hat S[\alpha])=H_{d-i} (\hat S [\beta])=0$. But $\alpha \subseteq V$ and $B=\hat S[V]$. So $\hat S [\alpha]=B[\alpha]$. Thus, $H_i (B[\alpha])=H_i (\hat S[\alpha])=0$. So, $\beta_i (B[\alpha];\mathbb{F})=0$ for all $\alpha \subseteq V(B)$. Dividing by $\binom{\# V}{\#\alpha}$ and adding over all $\alpha$, we get $\sigma_i(B)=0$. $\hfill{\Box}$

\begin{Corollary} If $S$ is a $k$-stacked $\mathbb{F}$-homology sphere of dimension $d \geq 2k+1$ then $\sigma_i (S; \mathbb{F})=0$ for $k\leq i \leq d-k-1$.
\end{Corollary}

{\bf Proof:} Take an $\mathbb{F}$-homology $k$-stacked $(d+1)$-ball $B$ such that $\text{Skel}_{d-k}(B)=\text{Skel}_{d-k}(S)$.
Then, for $i \leq d-k-1$ and any subset $\alpha$ of $V(S)=V(B)$, we have $H_i (S[\alpha])=H_i (B[\alpha])$, and therefore $\tilde \beta_i (S[\alpha]; \mathbb{F})=\tilde \beta_i (B[\alpha]; \mathbb{F})$. Dividing by $\binom{\# V(S)}{\# (\alpha)}$ and adding over all $\alpha$, we get $\sigma_i (S;\mathbb{F})=\sigma_i (B;\mathbb{F})$ for $i \leq d-k-1$. Therefore, this result follows from Lemma 1. $\hfill{\Box}$

Modifying Definition 2.1 in [2], we defined in [3] the mu-vector $(\mu_0,\ldots, \mu_d)$ (w.r.t. $\mathbb{F}$) of a $d$-dimensional simplicial complex $X$ as follows.
\begin{eqnarray*}
\mu_0 &=& \mu_0 (X;\mathbb{F})= \sum\limits_{x\in V(X)} \frac{1}{1+f_0 (\ell k (x,X))},\\
\mu_i &=& \mu_i (X; \mathbb{F}) = \sum\limits_{x\in V(X)} \frac{\delta_{i1}+\sigma_{i-1} (\ell k(x,X); \mathbb{F})}{1+f_0(\ell k(x,X)}, ~1 \leq i \leq d.
\end{eqnarray*}
Here, for $x\in V(X)$, the link $\ell k (x,X)$ of $x$ in $X$ is the subcomplex $\{\alpha : x \not\in \alpha, \alpha \sqcup \{x\} \in X\}$. Since the vertex links of any $k$-stacked $\mathbb{F}$-homology manifold of dimension $d+1$ with boundary (respectively of dimension $d$ without boundary) are $k$-stacked homology $d$-balls (respectively, $k$-stacked homology $(d-1)$-spheres), the following two results are immediate consequences of Lemma 1 and Corollary 2 (and the definition of the mu-vector).

\begin{Corollary} If $\Delta$ is a $k$-stacked $\mathbb{F}$-homology manifold with boundary then $\mu_i (\Delta; \mathbb{F})=0$ for $i \geq k+1$.
\end{Corollary}

\begin{Corollary} If $M$ is a $k$-stacked $\mathbb{F}$-homology manifold without boundary of dimension $d \geq 2k+2$ then $\mu_i (M; \mathbb{F})=0$ for $k+1 \leq i \leq d-k-1$.
\end{Corollary}

The importance of the mu-vector arises from the following two results from [3] (Theorems 1.6 and 1.7, op. cit.).

\begin{Theorem} If $M$ is an $\mathbb{F}$-homology $d$-manifold without boundary then $\mu_{d-i}(M;\mathbb{F})=\mu_i (M; \mathbb{F}), ~0 \leq i \leq d$.
\end{Theorem}

\begin{Theorem} Any simplicial complex $X$ of dimension $d$ satisifes:
\begin{enumerate}
\item[\rm(a)] $\sum\limits_{i=0}^\ell (-1)^{\ell-i} \mu_i (X; \mathbb{F}) \geq \sum\limits^\ell_{i=0} (-1)^{\ell-i} \beta_i (X; \mathbb{F})$ for $0 \leq \ell \leq d$.
\item[\rm(b)] Equality holds in (a) for some $\ell$ iff, for every induced subcomplex $Y$ of $X$, the $\mathbb{F}$-linear map $H_\ell (Y; \mathbb{F}) \rightarrow H_\ell (X; \mathbb{F})$, induced by the inclusion map, is injective.
\item[\rm(c)] $\mu_\ell (X; \mathbb{F}) \geq \beta_\ell (X;\mathbb{F})$ for $0 \leq \ell \leq d$.
\item[\rm(d)] Equality holds in (c) for some $\ell$ iff, for every induced subcomplex $Y$ of $X$, the $\mathbb{F}$-linear maps $H_{\ell-1}(Y; \mathbb{F}) \rightarrow H_{\ell-1} (X; \mathbb{F})$ and $H_\ell (Y; \mathbb{F}) \rightarrow H_\ell (X; \mathbb{F})$, induced by the inclusion map, are both injective.
\end{enumerate}
\end{Theorem}

Since, by Theorem 6, we have $0 \leq \beta_i (X; \mathbb{F}) \leq \mu_i (X;\mathbb{F})$, the following result due to Murai and Nevo [13, Theorem 4.4(i)] is immediate from Corollary 4 and Theorem 6 (c).

\begin{Theorem} (Murai and Nevo) Let $M$ be a $k$-stacked $\mathbb{F}$-homology manifold without boundary, of dimension $d\geq 2k+2$. Then $\beta_i (M; \mathbb{F})=0$ for $k+1 \leq i \leq d-k-1$.
\end{Theorem}

Similarly, the following result is immediate from Corollary 3 and Theorem 6 (c).

\begin{Theorem} If $\Delta$ is a $k$-stacked $\mathbb{F}$-homology manifold with boundary then $\beta_i (\Delta; \mathbb{F})=0$ for $i \geq k+1$.
\end{Theorem}

The following interesting consequence of Theorem 6 is also worth noting.

\begin{Corollary} Let $X$ be a simplicial complex of dimension $d$ and let $0 < k < d$. Suppose $\mu_{k-1} (X; \mathbb{F})=\beta_{k-1} (X; \mathbb{F})$ and $\mu_{k+1} (X; \mathbb{F}) =\beta_{k+1} (X; \mathbb{F})$. Then $\mu_k (X; \mathbb{F})=\beta_k (X;\mathbb{F})$.
\end{Corollary}

{\bf Proof:} Since $\mu_{k-1}=\beta_{k-1}$ and $\mu_{k+1} =\beta_{k+1}$, Theorem 6 (d) implies that, for every induced subcomplex $Y$ of $X$, the maps $H_{k-1} (Y; \mathbb{F}) \rightarrow H_{k-1} (X;\mathbb{F})$ and $H_k (Y; \mathbb{F}) \rightarrow H_k(X; \mathbb{F})$ are both injective. Theorem 6 (b) now implies that $\sum\limits_{i=0}^{k-1} (-1)^{k-1-i} \mu_i =\sum\limits^{k-1}_{i=0} (-1)^{k-1-i} \beta_i$ and $\sum\limits^k_{i=0} (-1)^{k-i} \mu_i =\sum\limits^k_{i=0} (-1)^{k-i} \beta_i$. Adding these two equations, we get $\mu_k =\beta_k$. $\hfill{\Box}$

The following result is also immediate from Theorem 6 (d) and the definition of $\mathbb{F}$-tightness.

\begin{Theorem} A $d$-dimensional simplicial complex $X$ is $\mathbb{F}$-tight iff $X$ is connected and $\mu_i (X;\mathbb{F})=\beta_i (X;\mathbb{F})$ for all $i$, $0 \leq i \leq d$.
\end{Theorem}

Recall that a simplicial complex $X$ is said to be $(k+1)$-neighbourly if any $k+1$ vertices of $X$ form a face of $X$ (of dimension $k$), i.e., if $f_k (X) =\binom{f_0(X)}{k+1}$. Since $k>0$, any $(k+1)$-neighbourly simplicial complex is connected. If $k>1$, then such a complex is simply connected and hence $\mathbb{F}$-orientable. The following result is an elementary consequence of the definitions of the mu- and sigma-vectors (and the obvious fact that if $X$ is $(k+1)$-neighbourly then all the vertex links $L$ of $X$ are $k$-neighbourly and hence so are all the induced subcomplexes of $L$) see [2, Lemma 3.9].

\begin{Lemma} If $X$ is a $(k+1)$-neighbourly simplicial complex then $\mu_i (X; \mathbb{F})=0 =\beta_i (X;\mathbb{F})$ for $1 \leq i \leq k-1$ and $\mu_0 (X; \mathbb{F})=1=\beta_0 (X; \mathbb{F})$.
\end{Lemma}

Now, the following result from [10] is immediate.

\begin{Theorem} (K\"{u}hnel) Let $M$ be a $(k+1)$-neighbourly $\mathbb{F}$-homology manifold without boundary, of dimension $2k$. Suppose $M$ is $\mathbb{F}$-orientable. Then $M$ is $\mathbb{F}$-tight.
\end{Theorem}

{\bf Proof:} By Lemma 11, we have $\mu_i =\beta_i$ for $0 \leq i \leq k-1$. By Theorem 5, $\mu_{2k-i}=\mu_i$, while (since $M$ is $\mathbb{F}$-orientable), Poincare duality gives $\beta_{2k-i} =\beta_i$. Therefore, we get $\mu_i =\beta_i$ for $k+1 \leq i \leq 2k$. So, $\mu_i =\beta_i$ for $i =k \pm 1$. Hence, by Corollary 9, $\mu_k =\beta_k$. Thus $\mu_i =\beta_i$ for $0 \leq i \leq 2k$. Hence by Theorem 10, $M$ is $\mathbb{F}$-tight. $\hfill{\Box}$

The following is one of the main results of this paper.

\begin{Theorem} Let $\Delta$ be a $(k+1)$-neighbourly $k$-stacked $\mathbb{F}$-homology manifold with boundary. Then $\Delta$ is $\mathbb{F}$-tight.
\end{Theorem}

{\bf Proof:} By Lemma 11, $\mu_i =\beta_i$ for $i \leq k -1$. Also by Corollary 3, $\mu_i =\beta_i$ for $i \geq k+1$. Hence by Corollary 9, $\mu_k =\beta_k$ as well. So, by Theorem 10, we are done. $\hfill{\Box}$

Our second main result is the following.

\begin{Theorem} Let $M$ be an $\mathbb{F}$-orientable $(k+1)$-neighbourly $k$-stacked $\mathbb{F}$-homology manifold without boundary, of dimension $d\neq 2k+1$. Then $M$ is $\mathbb{F}$-tight.
\end{Theorem}

{\bf Proof:} First suppose $d \leq 2k-1$. Then each vertex link of $M$ is a $k$-nieghbourly $\mathbb{F}$-homology sphere of dimension $\leq 2k-2$. Using the Dehn-Sommerville equations, it is easy to see that the only such homology spheres are the boundary complexes of simplices. It follows that $M$ itself is the boundary complex of a simplex. In this case $M$ is trivially $\mathbb{F}$-tight.

So, let $d \geq 2k$. If $d=2k$, then the result follows from Theorem 12. So, we may assume $d\geq 2k+2$. By Lemma 11, we have $\mu_i =\beta_i$ for $0 \leq i \leq k-1$. Hence, as in the proof of Theorem 12, duality yields $\mu_i =\beta_i$ for $d-k+1 \leq i \leq d$. Also, by Corollary 4, $\mu_i =\beta_i$ for $k+1 \leq i \leq d-k-1$. Thus, $\mu_i =\beta_i$ for all $i$, except possibly for $i=k$ and $i=d-k$. Now, Corollary 9 implies that $\mu_i =\beta_i$ for $i=k$, $d-k$ as well. Hence by Theorem 10, $M$ is $\mathbb{F}$-tight. $\hfill{\Box}$

Theorem 14 provides an affirmative answer to Question 5.5 in [8], except in the case $d=2k+1$.
Theorem 14 has the same statement as Theorem 3.11 (a) in [2] except that in the latter result, we had the hypothesis ``locally $k$-stellated" in place of ``$k$-stacked". But $k$-stellated spheres of dimension $\geq 2k-1$ are $k$-stacked by [1, Theorem 2.9], and locally $k$-stacked closed homology manifolds of dimension $\geq 2k+2$ are $k$-stacked by [1, Theorem 2.20] = [13, Theorem 4.6]. Thus Theorem 14 generalizes [2, Theorem 3.11 (a)]. We are unable to obtain a similar generalization of Theorem 3.11 (b) of [2]. But, undoubtedly, such a generalization should exist. Thus we posit:

\noindent {\bf Conjecture A:} Any $\mathbb{F}$-orientable $k$-stacked and $(k+1)$-neighbourly $\mathbb{F}$-homology manifold without boundary (of dimension $2k+1$, the only open case) is $\mathbb{F}$-tight.

Note that the reason for our inability to handle the case $d=2k+1$ is that Corollary 2 says nothing about the sigma vectors of $k$-stacked homology spheres of dimension $2k$. In [3, Conjecture 1] we made a sweeping conjecture on the sigma vectors of homology spheres. This conjecture includes:

\noindent{\bf Conjecture B:} If $S$ is a $k$-neighbourly $k$-stacked $\mathbb{F}$-homology sphere of dimension $2k$, then $\sigma_{k-1} (S; \mathbb{F})=\frac{\binom{m-k-2}{k+1}}{\binom{2k+3}{k+1}} -\delta_{k,1}$, where $m=f_0(S)$.

\begin{Theorem} Conjecture B $\Rightarrow$ Conjecture A.
\end{Theorem}

\noindent{\bf Sketch of proof:} Let $M$ be an $\mathbb{F}$-orientable $(2k+1)$-dimensional $k$-stacked and $(k+1)$-neighbourly $\mathbb{F}$-homology manifold without boundary. By Lemma 11 and duality, we have $\mu_i =\beta_i$ for $0 \leq i \leq k-1$, and for $k+2 \leq i \leq 2k+1$. Since, by duality, $\mu_{k+1} =\mu_k$ and $\beta_{k+1}=\beta_k$, it suffices to show that $\mu_k=\beta_k$. Then Theorem 10 would imply that $M$ is $\mathbb{F}$-tight , as claimed in Conjecture A.

By assumption, Conjecture B is valid of all the vertex links of $M$. Thus, for each $x$ in $V(M)$, $f_0(\ell k(x,M))=m-1$ and $\sigma_{k-1} (\ell k (x,M))=\frac{\binom{m-k-3}{k+1}}{\binom{2k+3}{k+1}}-\delta_{k,1}$ where $m=f_0(M)$. Therefore, we get $\mu_k =\frac{\binom{m-k-3}{k+1}}{\binom{2k+3}{k+1}}$. But Theorem 3.1 and Proposition 5.2 of [13] imply that for $M$ satisfying our hypotheses, $\beta_k$ is given by the same formula (see [3, Theorem 1.13] for the general statement). Thus $\mu_k =\beta_k$. $\hfill{\Box}$

The case $k=1$ of Conjecture B was proved by Burton et al in [6, Theorem 1.1]. Thus Conjecture A is true for $k=1$. This may also be deduced from [2, Theorem 3.11 (b)] and [1, Corollary 2.17].

\vskip 1em
\noindent{\bf Examples:} The then known examples of tight manifolds without boundary were surveyed by K\"{u}hnel and L\"{u}sz [11]. This list was updated in [2, Example 3.15]. Apart from the boundary complexes of simplices, this list includes the boundaries of the manifolds listed in (a) and (b) below. Excepting these series (and the 2-neighbourly 2-manifolds without boundary, which are tight by Theorem 12), only twenty nine sporadic examples of tight triangulated manifolds without boundary are known. Surprisingly, Theorems 13 and 14 notwithstanding, we do not know of any tight homology manifolds which are not triangulated (indeed combinatorial) manifolds. Perhaps these theorems will aid us in searching for such examples.

Apart from the trivial example of the standard $d$-ball (i.e., the face complex of the geometric $d$-simplex) only some low dimensional examples (cf. Banchoff [4, 5]) of tight triangulated manifolds with boundary appear to have been known prior to this paper. Note that Theorem 13 implies that any $k$-stacked $(k+1)$-neighbourly triangulated manifold with boundary is tight with respect to all fields. We have the following examples.

\begin{enumerate}
\item[(a)] This construction is due to K\"{u}hnel [9]. For $d\geq 2$, let ${\cal H}(d)$ be the pure simplicial complex of dimension $d$ whose facets are the $(d+1)$-vertex paths in a $(2d+1)$-vertex cyclic graph. Then ${\cal H}(d)$ is an 1-stacked 2-neighbourly triangulated handle body. Thus ${\cal H}(d)$ is tight (w.r.t. all fields) for $d\geq 2$. In particular, as Banchoff observed in [5], the 5-vertex M\"{o}bius band ${\cal H}(2)$ is tight.
\item[(b)] For each $d\geq 3$, Datta and Singh [7] constructed two non-isomorphic triangulated $d$-manifolds with boundary, named ${\cal M}_n^d$ and ${\cal N}_n^d$. Both are 1-stacked 2-neighbourly, with $d^2+3d+1$ vertices. They are tight with respect to all fields.
\item[(c)] The list in [2, Example 3.15] includes nineteen sporadic examples of triangulated manifolds without boundary which are $k$-stacked $(k+1)$-neighbourly for some $k \geq 1$. Each of them is the boundary of a tight manifold with boundary.
\end{enumerate}


\begin{thebibliography}{999}
\bibitem{} B. Bagchi and B. Datta, On $k$-stellated and $k$-stacked spheres, Discrete Math. 313 (2013), 2318-2329.
\bibitem{} B. Bagchi and B. Datta, On stellated spheres and a tightness criterion for combinatorial manifolds, Euro. J. Combin. 36 (2014), 294-313.
\bibitem{} B. Bagchi, The mu vector, Morse inequalities, and a generalized lower bound theorem for locally tame combinatorial manifolds, Preprint, arXiv : 1405.5675.
\bibitem{} T.F. Banchoff, Tightly embedded 2-dimensional polyhedral manifolds, Amer. J. Math. 87 (1965), 462-472.
\bibitem{} T.F. Banchoff, Tight polyhedral Klein bottles, projective planes and M\"{o}bius bands, Math. Ann. 207 (1974), 233-243.
\bibitem{} B.A. Burton, B. Datta, N. Singh and J. Spreer, Separation index of graphs and stacked 2-spheres, Preprint, arXiv : 1403.5862.
\bibitem{} B. Datta and N. Singh, An infinite family of tight triangulations of manifolds, J. Combin. Theory (A) 120 (2013), 2148-2163.
\bibitem{} F. Effenberger, Stacked polytopes and tight triangulations of manifolds, J. Combin. Theory (A) 118 (2011), 1843-1862.
\bibitem{} W. K\"{u}hnel, Higher dimensional analogues of Csaszar's torus, Results in Mathematics 9 (1986), 95-106.
\bibitem{} W. K\"{u}hnel, Tight polyhedral manifolds and tight triangulations, Lecture notes in Mathematics 1612, Springer Verlag, Berlin, 1995.
\bibitem{} W. K\"{u}hnel and F.H. L\"{u}sz, A census of tight triangulations, Period Math. Hungar. 39 (1999), 161-183.
\bibitem{} S. Murai and E. Nevo, On the generalized lower bound conjecture for polytopes and spheres, Acta Math. 210 (2013), 185-202.
\bibitem{} S. Murai and E. Nevo, On $r$-stacked triangulated manifolds, J. Algebraic Combin. 39 (2014), 373-388.
\end{thebibliography}
\end{document}